\documentclass[12pt]{amsart}

\usepackage{latexsym,amssymb,amsthm,upref,mathrsfs}
\usepackage{amsmath}
\usepackage{cite}
\usepackage[top=0.8in, left=1in, right=1in, bottom=1in]{geometry}
\usepackage{listings}
\newcommand\Fontv{\fontsize{8}{8.2}\selectfont}

\lstdefinelanguage{GAP}{%
  morekeywords={%
    Assert,Info,IsBound,QUIT,%
    TryNextMethod,Unbind,and,break,%
    continue,do,elif,%
    else,end,false,fi,for,%
    function,if,in,local,global,%
    mod,not,od,or,%
    quit,rec,repeat,return,%
    then,true,until,while,proc,options%
  },%
  sensitive,%
  morecomment=[l]\#,%
  morestring=[b]",%
  morestring=[b]',%
}[keywords,comments,strings]

\newcommand{\cl}{C \kern -0.1em \ell}
\DeclareMathOperator{\hotimes}{\Hat{\otimes}}

\makeatletter
\newcommand\GL{\mathop{\operator@font GL}\nolimits}
\makeatother

\makeatletter
\newcommand\Mat{\mathop{\operator@font Mat}\nolimits}
\makeatother

\makeatletter
\newcommand\sgn{\mathop{\operator@font sgn}\nolimits}
\makeatother

\makeatletter
\newcommand\spn{\mathop{\operator@font span}\nolimits}
\makeatother

\makeatletter
\newcommand\tr{\mathop{\operator@font tr}\nolimits}
\makeatother

\newcommand{\BR}{\mathbb{R}}
\newcommand{\BC}{\mathbb{C}}
\newcommand{\BH}{\mathbb{H}}

\newcommand{\BZ}{\mathbb{Z}}
\newcommand{\BF}{\mathbb{F}}

\newcommand{\BK}{\mathbb{K}}

\newcommand{\Pin}{\mathbf{Pin}}
\newcommand{\Lipschitz}{\mathbf{\Gamma}}

\newcommand{\be}{\mathbf{e}}
\newcommand{\bi}{\mathbf{i}}
\newcommand{\bj}{\mathbf{j}}

\newcommand{\bx}{\mathbf{x}}

\newcommand{\ed}{\end{document}}

\def\ve{\varepsilon}
\def\vp{\varphi}

\newcommand{\ta}[2]{#1_#2\tilde{\phantom{.}}}

\newcommand{\tp}{\ta{T}{\ve}}
\newcommand{\cb}[1]{\mathcal{#1}}

\newcommand{\Gpq}[2]{G_{#1,#2}}
\newcommand{\Gpqf}[3]{G_{#1,#2}(#3)}
\newcommand{\clpq}[2]{\cl_{#1,#2}}
\newcommand{\clpqr}[3]{\cl_{#1,#2,#3}}

\newcommand{\Kpqf}[3]{K_{#1,#2}(#3)}
\newcommand{\Tpqf}[3]{T_{#1,#2}(#3)}

\newcommand{\fpower}[1]{{}^2 \kern -0.2em #1} 

\newcommand{\iu}{{\underline{i}}}
\newcommand{\ju}{{\underline{j}}}

\newcommand{\ua}{\underline{a}}
\newcommand{\ub}{\underline{b}}
\newcommand{\uk}{\underline{k}}

\newcommand{\End}{\mathrm{End}}

\theoremstyle{plain}
\newtheorem{theorem}{Theorem}
\newtheorem{corollary}{Corollary}
\newtheorem{lemma}{Lemma}
\newtheorem{proposition}{Proposition}

\newtheorem*{StructureTheorem}{Structure Theorem}
\theoremstyle{definition}
\newtheorem{definition}{Definition}
\newtheorem{example}{Example}
\newtheorem{remark}{Remark}

\begin{document}
\vspace*{-4ex}
\title{On Clifford Algebras and Related\\ Finite Groups and Group Algebras}

\maketitle
\begin{center}
\author{{\bf Rafa\l \ Ab\l amowicz}\,$^a$}\\
\vspace{15pt}
\small
\thispagestyle{empty}
$^a$ Department of Mathematics, Tennessee Technological University \\
Cookeville, TN 38505, U.S.A. \\
\texttt{rablamowicz@tntech.edu}, \texttt{http://math.tntech.edu/rafal/}
\end{center}

\noindent
\begin{abstract}
Albuquerque and Majid~\cite{albuquerque} have shown how to view Clifford algebras $\cl_{p,q}$ as twisted group rings whereas Chernov has observed~\cite{chernov} that Clifford algebras can be viewed as images of group algebras of certain $2$-groups modulo an ideal generated by a nontrivial central idempotent. Ab\l amowicz and Fauser~\cite{ablamowicz1,ablamowicz2,ablamowicz3} have introduced a special transposition anti-automorphism of $\cl_{p,q}$, which they called a ``transposition", which reduces to reversion in algebras $\cl_{p,0}$ and to conjugation in algebras $\cl_{0,q}$. The purpose of this paper is to bring these concepts together in an attempt to investigate how the algebraic properties of real Clifford algebras, including their periodicity of eight, are a direct consequence of the central product structure of Salingaros vee groups viewed as $2$-groups.
\end{abstract}

\noindent
{\small

}
\vspace*{-4ex}
\tableofcontents

\section{Introduction}
The main goal of this survey paper is to show how certain finite groups, in particular, Salingaros vee groups~\cite{salingaros1,salingaros2,salingaros3}, and elementary abelian group 
$(\BZ_2)^n = \BZ_2 \times \cdots \times \BZ_2$ ($n$-times), and their group algebras and twisted groups algebras, arise in the context of Clifford algebras $\clpq{p}{q}.$ 

Chernov's observation~\cite{chernov} that Clifford algebras $\clpq{p}{q}$ can be viewed as images of (non-twisted) group algebras of suitable $2$-groups, conjectured to be Salingaros vee groups~\cite{walley},  allows one to gain a new viewpoint on these algebras and to relate classical group-theoretical results\cite{dornhoff,gorenstein,mckay}, in particular, on finite $2$-groups, to the theory of Clifford algebras. Salingaros classified the groups $\Gpq{p}{q}$ -- referred to as \textit{Salingaros vee groups} -- into five non-isomorphic classes $N_{2k-1},$ $N_{2k},$ $\Omega_{2k-1},$ $\Omega_{2k},$ and $S_{k}$. These groups, according to the theory of finite $2$-groups~\cite{dornhoff,mckay}, are central products of extra special groups $D_8$ -- the dihedral group and $Q_8$ -- the quaternionic group, both of order~$8$, and their centers $\BZ_2$, $\BZ_2 \times \BZ_2$, or $\BZ_4$. Thus, the properties of these groups and the fact that they fall into the five classes, is reflected by the fact that Clifford algebras $\clpq{p}{q}$ also fall into five isomorphism classes which is well known~\cite{chevalley,lam,lounesto} and references therein. The structure theorem on these algebras is recalled in Appendix~A. Furthermore, the ``periodicity of eight" of Clifford algebras viewed as the images of Salingaros vee groups, seems to be related to, if not predicted by, the structure of these groups and their group algebras. Thus, Section~2 is devoted to this approach to Clifford algebras.

Section~3 is devoted to a review of the basic properties of Salingaros vee groups $\Gpq{p}{q}$ appearing as finite subgroups of the group of units $\clpq{p}{q}^{\times}.$ Furthermore, we will review certain important subgroups of these groups appearing in the context of certain stabilizer groups of primitive idempotents in $\clpq{p}{q} $~\cite{ablamowicz2,ablamowicz3}.

Section~4 is devoted to a review of the central product structure of Salingaros vee groups.

In Section~5, we recall how the elementary abelian group $(\BZ_2)^n$ appears in the context of defining Clifford product on the set of monomials $\be_{\ua}$ indexed by binary $n$-tuples~$\ua$ from $(\BZ_2)^n$. In this first context, Walsh functions -- essentially, irreducible characters of $(\BZ_2)^n$ -- and Gray code -- as a certain isomorphism of $(\BZ_2)^n$ -- are used to define the $\clpq{p}{q}$ algebra product~\cite[Page 284]{lounesto} and references therein. In particular, a formula given by Lounesto dates back to 1935 and is being attributed to Brauer and Weyl~\cite{brauer}. It will be shown how this formula, applicable only to real Clifford algebras $\clpq{p}{q}$ over quadratic vector spaces $(V,Q)$ with a non-degenerate quadratic form $Q$ of signature $(p,q)$, and for an orthonormal set of basis elements (group generators), can be easily extended to Clifford algebras $\clpqr{p}{q}{r}$ for degenerate quadratic form $Q$ with $\dim V^{\perp}=r.$ 

Finally, in Section~6, we briefly recall the group $(\BZ_2)^n$ as it appears again in the context of the Clifford algebra $\clpq{p}{q}$ as a twisted group algebra $\BR^t[(\BZ_2)^n]$ viewed as a Hopf algebra with a certain quasi-triangular structure~\cite{albuquerque,downs}. This structure is needed to twist the commutative product in the group algebra $\BR[(\BZ_2)^n]$ in a manner similar to the Brauer and Weyl formula, so that the twisted product is the Clifford product in $\clpq{p}{q}.$ It is recalled that the ``transposition" anti-involution of $\clpq{p}{q}$ introduced in~\cite{ablamowicz1,ablamowicz2,ablamowicz3} is actually the antipode in the Hopf algebra $\BR^t[(\BZ_2)^n]$.\footnote{We remark that twisted group rings can also be described as certain special Ore extensions known as skew polynomial rings~\cite{bueso}.}

Our standard references on the group theory are~\cite{dornhoff, gorenstein, rotman}; in particular, for the theory of $p$-groups we rely on~\cite{mckay}; for Clifford algebras we use~\cite{chevalley,lam,lounesto} and references therein; on representation theory we refer to~\cite{james}; and for the theory of Hopf algebras we refer to~\cite{majid}.  

\section{Clifford Algebras as Images of Group Algebras}
Using Chernov's idea \cite{chernov}, in this section we want to show how Clifford algebras $\clpq{p}{q}$ can be viewed as images of group algebras $\BR[G]$ of certain $2$-groups. It is conjectured~\cite{walley} that the group~$G$, up to an isomorphism, is the Salingaros vee group 
$\Gpq{p}{q}$~\cite{salingaros1,salingaros2,salingaros3}. These groups, and their subgroups, have been recently discussed in \cite{ablamowicz2,ablamowicz3,brown,maduranga,maduranga2}.

\begin{definition}
Let $G$ be a finite group and let $\BF$ be a field\footnote{Usually, $\BF=\BR$ or $\BC$ although finite fields are also allowed. In this paper, we will be looking at the real Clifford algebras 
$\clpq{p}{q}$ as images of real group algebras or as real twisted group algebras.}. Then the \textit{group algebra} $\BF[G]$ is the vector space 
\begin{gather}
\BF[G] = \left\{\sum_{g \in G} \lambda_g g, \; \lambda_g \in \BF\right\}
\end{gather}
with multiplication defined as
\begin{gather}
\left(\sum_{g \in G} \lambda_g g\right)
\left(\sum_{h \in G} \mu_h h\right)=
\sum_{g,h \in G} \lambda_g \mu_h (gh)=
\sum_{g \in G} \sum_{h \in G} \lambda_h \mu_{h^{-1}g} g
\end{gather}
where all $\lambda_g,\mu_h \in \BF.$~\cite{james}
\end{definition}
Thus, group algebras are associative unital algebras with the group identity element playing the role of the algebra identity. In the theory of representations of finite groups, all irreducible inequivalent representations are related to a complete decomposition of the group algebra over $\BC$ viewed as a \textit{regular} $\BC$-module (cf. \cite[Maschke Theorem]{james}). The theory is rich on its own. The theory of group characters can then be derived from the representation 
theory~\cite{james}, or, as it is often done, from the combinatorial arguments and the theory of characters of the symmetric group~\cite{sagan}. Since in this survey we are only interested in finite groups, we just recall for completeness that every finite group is isomorphic to  a subgroup of a symmetric group~\cite{rotman}.

We begin by recalling a definition of a $p$-group.

\begin{definition}
Let $p$ be a prime. A group $G$ is a \textit{$p$-group} if every element in $G$ is of order $p^k$ for some $k \geq 1$.   
\end{definition}

Note that any finite group $G$ of order $p^n$ is a $p$-group. A classical result states that a center of any $p$-group is nontrivial, and, by Cauchy's theorem we know that every finite $p$-group has an element of order $p$. Thus, in particular, the center of any finite $p$-group has an element of 
order~$p$~\cite{dornhoff, gorenstein, rotman}. In the following, we will be working only with finite $2$-groups such as, for example, the group $(\BZ_2)^n$ and Salingaros vee groups $\Gpq{p}{q}$ of order $2^{1+p+q}.$

Two important groups in the theory of finite $2$-groups and hence in this paper, are the \textit{quaternionic group} $Q_8$ and the \textit{dihedral group} $D_8$ (the symmetry group of a square under rotations and reflections), both of order $|Q_8|=|D_8|=8.$ These groups have the following presentations:

\begin{definition}
The \textit{quaternionic group} $Q_8$ has the following two presentations:
\begin{subequations}
\begin{align}
Q_8 &= \langle a,b \mid a^4=1, a^2=b^2, bab^{-1}=a^{-1} \rangle \label{eq:q8a}\\
    &= \langle I,J,\tau \mid \tau^2=1, I^2=J^2=\tau, IJ=\tau J I\rangle \label{eq:q8b}
\end{align}
\end{subequations}
\end{definition}
\noindent
Thus, $Q_8=\{1,a,a^2,a^3,b,ab,a^2b,a^3b\}$ where the group elements have orders as follows: $|a^2|=2$, $|a|=|a^3|=|b|=|ab|=|a^2b|=|a^3b|=4,$  so the order structure of $Q_8$ is $(1,1,6),$\footnote{That is, $Q_8$ has one element of order 1; one element of order 2; and six elements of order 4.} and the center $Z(Q_8)=\{1,a^2\} \cong \BZ_2$. Here, we can choose $\tau=a^2.$ While the presentation~(\ref{eq:q8a}) uses only two generators, for convenience and future use, we prefer presentation~(\ref{eq:q8b}) which explicitly uses a central element $\tau$ of order~$2$.

\begin{definition}
The \textit{dihedral group} $D_8$ (the symmetry group of a square) has the following two presentations:
\begin{subequations}
\begin{align}
D_8 &= \langle a,b \mid a^4=b^2=1, bab^{-1}=a^{-1} \rangle \label{eq:d8a}\\
    &= \langle \sigma, \tau \mid \sigma^4=\tau^2=1, \tau \sigma \tau^{-1}=\sigma^{-1} \rangle  \label{eq:d8b}
\end{align}
\end{subequations}
\end{definition}
\noindent
Thus, $D_8=\{1,a,a^2,a^3,b,ab,a^2b,a^3b\}$ where $|a^2|=|b|=|ab|=|a^2b|=|a^3b|=2,$ $|a|=|a^3|=4,$ the order structure of $D_8$ is $(1,5,2),$ and $Z(D_8)=\{1,a^2\} \cong \BZ_2$. Here, we can choose $\tau=b,$ $\sigma=a$, hence, $\sigma^2 \in Z(D_8).$ That is, $\sigma^2$ is our central element of order~$2$, and our preferred presentation of $D_8$ is~(\ref{eq:d8b}).

In the following two examples, we show how one can construct the Clifford algebra 
$\clpq{0}{2} \cong \BH$ (resp.  $\clpq{1}{1}$) as an image of the group algebra of $Q_8$ (resp. $D_8$). 
\begin{example}(Constructing $\BH \cong \cl_{0,2}$ as $\BR[Q_8]/\mathcal{J}$)\\
Define an algebra map $\psi$ from the group algebra $\BR[Q_8] \rightarrow \BH = \spn_\BR\{1,\bi,\bj,\bi\bj \}$ as follows:
\begin{gather}
1 \mapsto 1,\quad \tau \mapsto -1,\quad I \mapsto \bi,\quad J  \mapsto \bj, 
\end{gather}
Then, $\cb{J}=\ker \psi = (1+\tau)$ for the central element $\tau$ of order~$2$ in $Q_8$\footnote{Here, $(1+\tau)$ denotes an ideal in $\BR[Q_8]$ generated by $1+\tau$. Note that the two elements $\frac12(1\pm\tau)$ are idempotents which provide an \textit{orthogonal decomposition} of the unity in $\BR[Q_8]$.}, so $\dim_{\BR} \cb{J} = 4$ and $\psi$ is surjective. Let $\pi:\BR[Q_8]\rightarrow \BR[Q_8]/\cb{J}$ be the natural map $u \mapsto u+\cb{J}.$ There exists an isomorphism 
$\vp: \BR[Q_8]/\cb{J} \rightarrow \BH$ such that $\vp \circ \pi = \psi$ and 
\begin{gather*}
\pi(I^2) = I^2+\cb{J} = \tau+\cb{J} \mbox{ and } \vp(\pi(I^2)) = \psi(\tau)=-1
=(\psi(I))^2=\bi^2,\\
\pi(J^2) = J^2+\cb{J} = \tau+\cb{J} \mbox{ and } \vp(\pi(J^2)) = \psi(\tau)=-1
=(\psi(J))^2=\bj^2,\\
\pi(IJ+JI)=IJ+JI+\cb{J}=(1+\tau)JI+\cb{J} = \cb{J}  \mbox{ and }\\ 
\vp(\pi(IJ+JI)) = \psi(0)=0 = \psi(I)\psi(J)+\psi(J)\psi(I)=\bi\bj+\bj\bi.
\end{gather*}
Thus, $\BR[Q_8]/\cb{J} \cong \psi(\BR[Q_8]) = \BH \cong \cl_{0,2}$ provided the central element $\tau$ is mapped to $-1$ (see also~\cite{chernov}).
\label{ex:example1}
\end{example}

\begin{example} (Constructing $\cl_{1,1}$ as $\BR[D_8]/\mathcal{J}$)\\
Define an algebra map $\psi$ from the group algebra $\BR[D_8] \rightarrow \cl_{1,1}$ such that:
\begin{gather}
1 \mapsto 1,\quad \tau \mapsto \be_1,\quad \sigma \mapsto \be_2,
\label{eq:rq8} 
\end{gather}
where $\cl_{1,1}=\spn_\BR \{1,\be_1,\be_2,\be_1\be_2 \}$. Then, $\ker \psi = (1+\sigma^2)$ where $\sigma^2$ is a central element of order~$2$ in $D_8$. Let $\cb{J} = (1+\sigma^2)$. Thus, $\dim_{\BR} \cb{J} = 4$ and 
$\psi$ is surjective. Let $\pi:\BR[D_8]\rightarrow \BR[D_8]/\cb{J}$ be the natural map 
$u \mapsto u+\cb{J}.$ There exists an isomorphism $\vp: \BR[D_8]/\cb{J} \rightarrow \cl_{1,1}$ such that $\vp \circ \pi = \psi$ and 
\begin{gather*}
\pi(\tau^2) = \tau^2+\cb{J} = 1+\cb{J} \mbox{ and } \vp(\pi(\tau^2)) = \psi(1)=1=\psi(\tau^2)=(\be_1)^2,\\
\pi(\sigma^2) = \sigma^2+\cb{J} \mbox{ and } \vp(\pi(\sigma^2)) = \psi(\sigma^2)=-1=
(\be_2)^2,\\
\pi(\tau\sigma+\sigma\tau)=\tau\sigma+\sigma\tau+\cb{J}=\sigma\tau(1+\sigma^2)+\cb{J}=\cb{J} \mbox{ and }\\ \vp(\pi(\tau\sigma+\sigma\tau) = \psi(0)=0=\psi(\tau)\psi(\sigma)+\psi(\sigma)\psi(\tau)=\be_1\be_2+\be_2\be_1.
\end{gather*}
Thus, $\BR[D_8]/\cb{J} \cong \cl_{1,1}$ provided the central element $\sigma^2$ is mapped to $-1.$
\label{ex:example2}
\end{example}
It is not difficult to modify Example~2 and construct $\clpq{2}{0}$ as the quotient algebra 
$\BR[D_8]/\cb{J}$ by changing only the definition of the algebra map $\psi$ given in~(\ref{eq:rq8}) to
\begin{gather}
1 \mapsto 1,\quad \tau \mapsto \be_1,\quad \sigma \mapsto \be_1\be_2,
\label{eq:rq88} 
\end{gather}
Then, the rest of Example~2 follows except that of course now $(\be_1)^2=(\be_2)^2=1$. Thus, one can construct $\cl_{2,0}$ as $\BR[D_8]/\mathcal{J}$ with again $\cb{J} = (1+\sigma^2).$

We remark that the fact that can use the group $D_8$ twice should not come as surprise since 
$\clpq{1}{1} \cong \clpq{2}{0}$ (as real Clifford algebras) due to one of the isomorphism theorems stating that 
$\clpq{p}{q} \cong \clpq{q+1}{p-1},$~\cite[Page 215]{lounesto} (see also~\cite{chevalley,lam,periodicity}) and that we only have, up to an isomorphism, two non-abelian groups of order eight, namely, $Q_8$ and $D_8.$ 

We summarize our two examples as follows. In preparation for Chernov's theorem~\cite{chernov}, notice 
that elements in each group $Q_8$ and $D_8$ can be written as follows:
\begin{itemize}
\item The quaternionic group $Q_8$:
$$
Q_8 = \{\tau^{\alpha_0}g_1^{\alpha_1}g_2^{\alpha_2} \mid \alpha_k \in \{0,1\},\, k=0,1,2\}
$$
where $\tau=a^2$ is the central element of order~$2$ in $Q_8$, $g_1=a$, and $g_2=b$. Thus, 
$$
(g_1)^2=a^2=\tau,  \quad (g_2)^2=b^2=a^2=\tau,  \quad \tau g_1g_2 =g_2g_1.
$$
Observe that $|g_1|=|g_2|=4$ and $\BR[Q_8]/\cb{J} \cong \cl_{0,2}$ where 
$\cb{J} = (1+\tau).$  
\item The dihedral group $D_8$:
$$
D_8 = \{\tau^{\alpha_0}g_1^{\alpha_1}g_2^{\alpha_2} \mid \alpha_k \in \{0,1\},\, k=0,1,2\}
$$
where $\tau=a^2$ is the central element of order~$2$ in $D_8$, $g_1=b$, and $g_2=a$. Thus, 
$$
(g_1)^2=b^2=1,  \quad (g_2)^2=a^2=\tau,  \quad \tau g_1g_2 =g_2g_1.
$$  
Observe that $|g_1|=2,$ $|g_2|=4$ and $\BR[D_8]/\cb{J} \cong \cl_{1,1}$ where 
$\cb{J} = (1+\tau).$
\end{itemize}

Chernov's theorem states the following.

\begin{theorem}[Chernov]
Let $G$ be a finite $2$-group of order $2^{1+n}$ generated by a central element~$\tau$ of order~$2$ and additional elements $g_1,\ldots,g_n,$ which satisfy the following relations:
\begin{subequations}
\begin{gather}
\tau^2=1, \quad (g_1)^2=\cdots=(g_p)^2=1, \quad (g_{p+1})^2=\cdots=(g_{p+q})^2=\tau,\label{eq:squares}\\
\tau g_j = g_j\tau, \quad g_ig_j=\tau g_jg_i, \quad i,j=1,\ldots,n=p+q, \label{eq:tau}
\end{gather}
\label{eq:both}
\end{subequations}
\hspace*{-1ex}so that $G=\{\tau^{\alpha_0}g_1^{\alpha_1}\cdots g_n^{\alpha_n} \mid \alpha_k \in \{0,1\},\, k=0,1,\ldots,n\}$. 
Let $\cb{J}=(1+\tau)$ be an ideal in the group algebra $\BR[G]$ and let $\cl_{p,q}$ be the universal real Clifford algebra generated by $\{\be_k\}, k=1,\ldots,n=p+q,$ where
\begin{subequations}
\begin{gather}
\be_i^2 = Q(\be_i) \cdot 1 = \ve_i \cdot 1 = 
\begin{cases}
1  & \textit{for $1 \leq i \leq p$;}\\
-1 & \textit{for $p+1 \leq i \leq p+q$;}
\end{cases}
          \label{eq:B1aa}\\
\be_i\be_j + \be_j \be_i = 0, \quad i \neq j,\quad  1 \le i,j \leq n. 
          \label{eq:B1bb}
\end{gather}
\label{eq:B11}
\end{subequations}
\noindent
\hspace*{-1ex}Then, (a) $\dim_{\BR}\cb{J}=2^n;$ (b) There exists a surjective algebra homomorphism $\psi$ from the group algebra $\BR[G]$ to $\cl_{p,q}$ so that $\ker \psi = \cb{J}$ and $\BR[G]/\cb{J} \cong \cl_{p,q}.$
\label{thm:thm2}
\end{theorem}

\begin{remark}
Chernov's theorem does not give the existence of the group $G$. It only states that should such group 
exist whose generators satisfy relations~(\ref{eq:both}), the result follows. It is not difficult to conjecture that the group~$G$ in that theorem is in fact the Salingaros vee group $\Gpq{p}{q}$, that is, $\BR[\Gpq{p}{q}]/\cb{J} \cong \clpq{p}{q}$ (see \cite{walley}). In fact, we have seen it in 
Examples~\ref{ex:example1} and \ref{ex:example2} above.
\end{remark}

\begin{proof}[Chernov's theorem]
Observe that $G=\{\tau^{\alpha_0}g_1^{\alpha_1}\cdots g_n^{\alpha_n}\} \mid
\alpha_k \in \{0,1\},\, k=0,1,\ldots,n\}$. The existence of a central element $\tau$ of order~$2$ is guaranteed by a well-known fact that the center of any $p$-group is nontrivial, and by Cauchy 
Theorem.~\cite{rotman} Define an algebra homomorphism $\psi:\BR[G] \rightarrow \cl_{p,q}$ such that
\begin{gather}
1 \mapsto 1, \quad \tau \mapsto -1, \quad g_j \mapsto \be_j, \quad j=1,\ldots,n.
\end{gather}
Clearly, $\cb{J} \subset \ker \psi$. Let $u \in \BR[G]$. Then,
\begin{gather}
u=\sum_{\alpha}\lambda_{\alpha}\tau^{\alpha_0}g_1^{\alpha_1}\cdots g_n^{\alpha_n}
 =u_{1}+\tau u_{2}
\end{gather} 
where
\begin{subequations} 
\begin{gather}
u_i = \sum_{\widetilde{\alpha}}\lambda_{\widetilde{\alpha}}^{(i)}g_1^{\alpha_1}\cdots g_n^{\alpha_n}, \quad i=1,2,\\
\alpha=(\alpha_0,\alpha_1,\ldots,\alpha_n) \in \BR^{n+1} \quad \mbox{and} \quad 
\widetilde{\alpha}=(\alpha_1,\ldots,\alpha_n) \in \BR^{n}.
\end{gather}
\end{subequations}
Thus, if $u \in \ker \psi$, then 
\begin{gather}
\psi(u) = 
\sum_{\widetilde{\alpha}}(\lambda_{\widetilde{\alpha}}^{(1)}-
                          \lambda_{\widetilde{\alpha}}^{(2)})
                          \be_1^{\alpha_1}\cdots\be_n^{\alpha_n} = 0
\end{gather}
implies $\lambda_{\widetilde{\alpha}}^{(1)}=\lambda_{\widetilde{\alpha}}^{(2)}$ since
$\{\be_1^{\alpha_1}\cdots\be_n^{\alpha_n}\}$ is a basis in $\cl_{p,q}$. Hence,
\begin{gather}
u=(1+\tau)\sum_{\widetilde{\alpha}}\lambda_{\widetilde{\alpha}}^{(1)}g_1^{\alpha_1}\cdots g_n^{\alpha_n} \in \cb{J}.
\end{gather}
Thus, $\dim_{\BR} \ker \psi=2^n,$ $\ker \psi = \cb{J}$, $\dim_{\BR}\BR[G]/\cb{J} = 2^{1+n}-2^n=2^n,$ so $\psi$ is surjective. Let $\vp:\BR[G]/\cb{J}\rightarrow \cl_{p,q}$ be such that $\vp \circ \pi = \psi$ where $\pi:\BR[G]\rightarrow \BR[G]/\cb{J}$ is the natural map. Then, since $\psi(g_j)=\be_j,$ $\pi(g_j)=g_j+\cb{J},$ we have 
$\vp(\pi(g_j)) = \vp(g_j+\cb{J})=\psi(g_j)=\be_j$ and
\begin{align}
\pi(g_j)\pi(g_i)+\pi(g_i)\pi(g_j)
&=(g_j+\cb{J})(g_i+\cb{J}) + (g_j+\cb{J})(g_i+\cb{J}) \notag \\
&=(g_jg_i+g_ig_j) + \cb{J} = (1+\tau)g_jg_i + \cb{J} = \cb{J}
\end{align}
because $g_ig_j=\tau g_jg_i$ in $\BR[G]$, $\tau$ is central, and $\cb{J}=(1+\tau).$ Thus, $g_j+\cb{J},g_i+\cb{J}$ anticommute in $\BR[G]/\cb{J}$ when $i \neq j.$ Also,
\begin{gather}
\pi(g_i)\pi(g_i) = (g_i+\cb{J})(g_i+\cb{J}) = (g_i)^2+\cb{J} =
\begin{cases}
1 + \cb{J},     & \text{$1 \leq i \leq p$;}\\ 
\tau + \cb{J},  & \text{$p+1 \leq i \leq n$;}
\end{cases}
\end{gather}
due to the relations~(\ref{eq:squares}) on $g_i$ in $G$. Observe, that 
\begin{gather}
\tau +\cb{J} = (-1)+(1+\tau) +\cb{J} = (-1)+\cb{J} \mbox{ in } \BR[G]/\cb{J}.
\end{gather} 
To summarize, the factor algebra $\BR[G]/\cb{J}$ is generated by the cosets $g_i+\cb{J}$ which satisfy these relations:
\begin{subequations}
\begin{gather}
(g_j+\cb{J})(g_i+\cb{J}) + (g_j+\cb{J})(g_i+\cb{J}) = \cb{J},\\ 
(g_i)^2+\cb{J} =
\begin{cases}
1 + \cb{J},     & \text{$1 \leq i \leq p$;}\\ 
(-1) + \cb{J},  & \text{$p+1 \leq i \leq n$;}
\end{cases}
\end{gather}
\end{subequations}
Thus, the factor algebra $\BR[G]/\cb{J}$ is a Clifford algebra isomorphic to $\cl_{p,q}$ provided 
$\cb{J}=(1+\tau)$ for the central element $\tau$ of order~$2$ in $G$.
\end{proof}

\section{Salingaros Vee Groups $G_{p,q} \subset \cl_{p,q}^{\times}$}

Let $\dim_\BR V=n$ and $Q$ be a non-degenerate quadratic form on~$V$:
\begin{align}
Q(\bx) &= \ve_1x_1^2 + \ve_2x_2^2 + \cdots + \ve_nx_n^2,
\label{eq:Q}     
\end{align}
$\ve_i = \pm 1$ and $\bx = x_1\be_1 + \cdots + x_n\be_n \in V$ for an orthonormal basis $\mathcal{B}_1 = \{\be_i,1\le i \le n\}$. $Q$ has an arbitrary signature $-n \le p-q \le n$ where $p$ (resp. $q$) denotes the number of $+1$'s (resp. $-1$'s) in~(\ref{eq:Q}), and $p+q=n$. Let $\cl_{p,q}$ be the universal Clifford algebra of $(V,Q)$ obtained, for example, via Chevalley's construction~\cite[Chapter 22]{lounesto}. 

Then, let $\cb{B}$ be the canonical basis of $\bigwedge V$ generated by $\cb{B}_1,$ 
$[n]=\{1,2,\ldots,n\}$ and denote arbitrary, canonically ordered subsets of $[n]$, by underlined Roman characters. The basis elements of $\bigwedge V$, or, of $\clpq{p}{q}$ due to the linear space isomorphism $\bigwedge V \rightarrow \clpq{p}{q}$, can be indexed by these finite ordered subsets as $\be_\iu = \wedge_{i \in \iu}\, \be_i$. 

Now, let $G_{p,q}$ be a finite group in any real Clifford algebra $\cl_{p,q}$ (simple or semisimple) with a with a binary operation being just the Clifford product, namely:
\begin{equation}
G_{p,q} = \{ \pm \be_\iu  \; | \; \be_\iu \in \cb{B} \; \mbox{with} \; 
\be_\iu \be_\ju \; \mbox{denoting the Clifford product}\}.
\end{equation}
Thus, $G_{p,q}$ may be presented as follows:
\begin{gather}
G_{p,q}=\langle -1,\be_1,\ldots,\be_n \mid \be_i\be_j = -\be_j\be_i\, \mbox{ for }\, i \neq j \, \mbox{ and }
\be_i^2= \pm 1
\rangle ,
\label{eq:Gpq}
\end{gather} 
where $\be_i^2=1$ for $1 \leq i \leq p$ and $\be_i^2=-1$ for $p+1 \leq i \leq n=p+q$. In the following, the elements $\be_{\underline{i}} = \be_{i_1}\be_{i_2}\cdots \be_{i_k}$ will be denoted for short as $\be_{i_1i_2\cdots i_k}$ for $k \ge 1$ while $\be_{\emptyset}$ will be denoted as $1,$ the identity element of $G_{p,q}$ (and $\cl_{p,q}).$ 

This $2$-group of order $2 \cdot 2^{p+q} = 2^{n+1}$ is known as \textit{Salingaros vee group} and it has been discussed, for example, by Salingaros~\cite{salingaros1,salingaros2,salingaros3}, Varlamov~\cite{varlamov1,varlamov}, Helmstetter~\cite{helmstetter1}, Ab\l amowicz and Fauser~\cite{ablamowicz2,ablamowicz3}, 
Maduranga and Ab\l amowicz~\cite{maduranga2}, and most recently by Brown~\cite{brown}. We should recall here that $G_{p,q}$ is a discrete subgroup of $\Pin(p,q) \subset \Lipschitz_{p,q}$ (Lipschitz group)  (Lounesto~\cite{lounesto}).  

In preparation for discussing properties of the groups $\Gpq{p}{q}$ and related to them subgroups, we recall a definition of the derived subgroup $G' \subset G$ and a proposition that gives some of its properties~\cite{rotman}.

\begin{definition}
If $G$ is a group and $x,y \in G,$ then their \textit{commutator} $[x,y]$ is the element $xyx^{-1}y^{-1}.$ If $X$ and $Y$ are subgroups of $G$, then the \textit{commutator subgroup} $[X,Y]$ of $G$ is defined by 
\begin{gather}
[X,Y]=\langle [x,y] \mid x \in X, y \in Y \rangle,
\end{gather}
that is, the group $[X,Y]$ is generated by all the commutators $[x,y].$ In particular, the 
\textit{derived subgroup} $G'$ of $G$ is defined as $G'=[G,G].$
\end{definition}
\begin{proposition}
Let $G$ be a group.
\begin{itemize}
\item[(i)] $G'$ is a normal subgroup of $G$, and $G/G'$ is abelian.
\item[(ii)] If $H$ is a normal subgroup of G and $G/H$ is abelian, then $G'\subseteq H.$
\end{itemize}
\end{proposition}

\subsection{Transposition Anti-Involution in $\clpq{p}{q}$}
\label{sub:sectt}

Let us now recall a definition and some of its basic properties of a special anti-involution $\tp$ in a Clifford algebra $\clpq{p}{q}$ referred to as ``transposition". This anti-involution was introduced in  \cite{ablamowicz1,ablamowicz2,ablamowicz3} where its properties were investigated at length. In particular, it allowed for an introduction of a reciprocal basis in a Clifford algebra $\clpq{p}{q}$ and, subsequently, a new spinor product on spinor spaces, and a classification of its (infinite) groups of invariance. In the following, we limit ourselves only to reviewing certain finite groups appearing in this context.   

\begin{definition}
The \textit{transposition} $\tp$  of $\cl_{p,q}$ is defined as:
\begin{gather}
\tp :\cl_{p,q} \rightarrow \cl_{p,q},\qquad 
    \sum_{\iu \in 2^{[n]}} u_\iu \be_{\iu} \mapsto
    \sum_{\iu \in 2^{[n]}} u_\iu (\be_{\iu})^{-1}
\label{eq:tpdef}
\end{gather}
\end{definition}
\noindent
It is the \textit{antipode} map $S$ known from the theory of group algebras $\BF[G]$
\begin{gather}
\BF[G] \rightarrow \BF[G], \qquad \sum_{g\in G} \lambda_gg \mapsto \sum_{g\in G} \lambda_g g^{-1}
\label{eq:S}
\end{gather}
viewed as Hopf algebras~\cite{majid}. Here are a few of its properties and a few finite related groups. For more details, see \cite{ablamowicz1,ablamowicz2,ablamowicz3}.
\begin{itemize}
\item $\tp$ is an anti-involution of $\cl_{p,q}$ which reduces to reversion in 
$\cl_{p,0}$ and to conjugation in $\cl_{0,q}$.
\item Depending on the value of $(p-q) \bmod 8$, where $(p,q)$ is the signature of $Q$, $\tp$ gives rise to transposition, Hermitian complex, and Hermitian quaternionic conjugation of spinor representation matrices. 
\item $\ta{T}{\ve}(\be_\iu) = \be_\iu^{-1}$ hence $\tp(\be_\iu) = \be_\iu$ (resp. $\tp(\be_\iu) = -\be_\iu$) when $(\be_\iu)^2=1$ (resp. $(\be_\iu)^2=-1$) (elements of order 2 and 4, respectively, in $G_{p,q}$).
\item $\tp(f) = f$ for any primitive idempotent $f$.
\item Let $S=\cl_{p,q}f$ be a spinor (minimal left) ideal in a simple algebra $\cl_{p,q}$ generated by a primitive idempotent~$f$. Then, $\tp$ defines a dual spinor space 
$S^\ast=\tp(S)$ and a $\BK$-valued, where $\BK=f\clpq{p}{q}f,$ spinor norm $(\psi,\phi) = \tp(\psi)\phi$ on $S$ invariant under (infinite) group $G_{p,q}^\ve$ (with $G_{p,q} < G_{p,q}^\ve$) different, in general, from spinor norms related to reversion and conjugation in $\cl_{p,q}$.
\item $G_{p,q}$ act transitively on a complete set $\cb{F}$, $|\cb{F}|=2^{q-r_{q-p}}$,  of mutually annihilating primitive idempotents where $r_i$ is the Radon-Hurwitz number. See a footnote in Appendix~A for a definition of $r_i$.
\item The normal stabilizer subgroup $G_{p,q}(f) \lhd G_{p,q}$ of $f$ is of order 
$2^{1+p+r_{q-p}}$ and monomials $m_i$ in its (non-canonical) left transversal together with $f$ determine a spinor basis in $S$.
\item The stabilizer groups $G_{p,q}(f)$ and the invariance groups $G_{p,q}^\ve$ of the spinor norm have been classified according to the signature $(p,q)$ for $(p+q) \leq 9$ in simple and semisimple algebras $\cl_{p,q}$.
\item $G_{p,q}$ permutes the spinor basis elements modulo the commutator subgroup $G'_{p,q}$ by left multiplication.
\item The ring $\BK=f\cl_{p,q}f$ is $G_{p,q}$-invariant.
\end{itemize}

\subsection{Important Finite Subgroups of $\cl_{p,q}^{\times}$}
In this section, we summarize properties and definitions of some finite subgroups of the group of invertible elements $\cl_{p,q}^{\times}$ in the Clifford algebra $\clpq{p}{q}.$ These groups were defined in~\cite{ablamowicz1,ablamowicz2,ablamowicz3}.
\begin{itemize}
\item $G_{p,q}$ -- Salingaros vee group of order $|G_{p,q}|=2^{1+p+q}$,
\item $G'_{p,q} = \{1,-1\}$ -- the commutator subgroup of $G_{p,q}$,
\item Let $\cb{O}(f)$ be the orbit of $f$ under the conjugate action of $G_{p,q}$, and let $G_{p,q}(f)$ be the stabilizer of~$f$. Let 
\begin{gather}
N=|\cb{F}| = [G_{p,q}:G_{p,q}(f)]=|\cb{O}(f)|=|G_{p,q}|/|G_{p,q}(f)|=2\cdot 2^{p+q}/|G_{p,q}(f)|
\end{gather}
then $N=2^k$ (resp. $N=2^{k-1}$) for simple (resp. semisimple) $\cl_{p,q}$ where 
$k=q-r_{q-p}$ and $[G_{p,q}:G_{p,q}(f)]$ is the index of $G_{p,q}(f)$ in $G_{p,q}$. 
\item $G_{p,q}(f) \lhd G_{p,q}$ and  $|G_{p,q}(f)|=2^{1+p+r_{q-p}}$ 
(resp. $|G_{p,q}(f)|=2^{2+p+r_{q-p}}$) for simple (resp. semisimple) $\cl_{p,q}.$

\item The set of commuting monomials $\cb{T}= \{\be_{\iu_1},\ldots,\be_{\iu_k}\}$ (squaring to $1$) in the primitive idempotent
$
f = \frac12(1\pm\be_{\iu_1}) \cdots \frac12(1\pm\be_{\iu_k}) 
$ 
is point-wise stabilized by $G_{p,q}(f).$
\item $T_{p,q}(f) := \langle \pm 1, \cb{T}\rangle 
\cong \Gpq{p}{q}' \times \langle \be_{\iu_1},\ldots, \be_{\iu_k} \rangle 
\cong \Gpq{p}{q}' \times (\BZ_2)^k,$ the \textit{idempotent group} of $f$ with 
$|\Tpqf{p}{q}{f}|=2^{1+k}$,
\item $\Kpqf{p}{q}{f} =  \langle \pm 1, m \mid m \in \cb{K}\rangle < \Gpqf{p}{q}{f}$ -- the \textit{field group} of where $f$ is a primitive idempotent in $\cl_{p,q}$,
$\BK=f\cl_{p,q}f$, and $\cb{K}$ is a set of monomials (a transversal) in~$\cb{B}$ which span~$\BK$ as a real algebra. Thus,
\begin{gather}
|\Kpqf{p}{q}{f}| = \begin{cases} 2, & p - q = 0,1,2 \bmod 8;\\
                                 4, & p - q = 3,7 \bmod 8;\\ 
                                 8, & p - q = 4,5,6 \bmod 8.
                   \end{cases}
\label{eq:orderKpqf}
\end{gather}
\item $G_{p,q}^{\ve} = \{g \in \cl_{p,q} \mid \tp(g)g=1\} $ (infinite group)
\end{itemize}

Before we state the main theorem from~\cite{ablamowicz3} that relates the above finite groups to the Salingaros vee groups, we recall the definition of a \textit{transversal}.

\begin{definition}
Let $K$ be a subgroup of a group $G$. A \textit{transversal} $\ell$ of $K$ in~$G$ is a subset of~$G$ consisting of exactly one element $\ell(bK)$ from every (left) coset $bK$, and with $\ell(K)=1$. 
\end{definition}

\begin{theorem}[Main Theorem]
Let $f$ be a primitive idempotent in $\cl_{p,q}$ and let $\Gpq{p}{q}$, $\Gpqf{p}{q}{f}$, $\Tpqf{p}{q}{f}$, $\Kpqf{p}{q}{f}$, and $\Gpq{p}{q}'$ be the groups defined above. Let $S=\cl_{p,q}f$ and $\BK=f\cl_{p,q}f$.
\begin{itemize}
\item[(i)] Elements of $\Tpqf{p}{q}{f}$ and $\Kpqf{p}{q}{f}$ commute.
\item[(ii)] $\Tpqf{p}{q}{f} \cap \Kpqf{p}{q}{f} = \Gpq{p}{q}' = \{\pm 1 \}$.
\item[(iii)] $\Gpqf{p}{q}{f} = \Tpqf{p}{q}{f}\Kpqf{p}{q}{f} = 
\Kpqf{p}{q}{f}\Tpqf{p}{q}{f}$.
\item[(iv)]  $|\Gpqf{p}{q}{f}| = |\Tpqf{p}{q}{f}\Kpqf{p}{q}{f}| =
\frac12 |\Tpqf{p}{q}{f}||\Kpqf{p}{q}{f}|$. 
\item[(v)] $\Gpqf{p}{q}{f} \lhd \Gpq{p}{q}$, $\Tpqf{p}{q}{f} \lhd \Gpq{p}{q}$,
and $\Kpqf{p}{q}{f} \lhd \Gpq{p}{q}$. In particular, $\Tpqf{p}{q}{f}$ and
$\Kpqf{p}{q}{f}$ are normal subgroups of $\Gpqf{p}{q}{f}$.
\item[(vi)] We have:
\begin{align}
  \Gpqf{p}{q}{f} /\Kpqf{p}{q}{f} &\cong \Tpqf{p}{q}{f} /\Gpq{p}{q}',\\
  \Gpqf{p}{q}{f} /\Tpqf{p}{q}{f} &\cong \Kpqf{p}{q}{f} /\Gpq{p}{q}'.
  \label{eq:conj6}
\end{align}
\item[(vii)] We have:
\begin{gather}
  (\Gpqf{p}{q}{f}/\Gpq{p}{q}')/(\Tpqf{p}{q}{f}/\Gpq{p}{q}')
   \cong \Gpqf{p}{q}{f}/\Tpqf{p}{q}{f} \cong \Kpqf{p}{q}{f}/\{\pm 1 \} 
  \label{eq:conj7}
\end{gather}
and the transversal of $\Tpqf{p}{q}{f}$ in $\Gpqf{p}{q}{f}$ spans $\BK$
over $\BR$ modulo~$f$.
\item[(viii)] The transversal of $\Gpqf{p}{q}{f}$ in $\Gpq{p}{q}$ spans $S$
over $\BK$ modulo~$f$.  
\item[(ix)] We have $(\Gpqf{p}{q}{f}/\Tpqf{p}{q}{f}) \lhd
(\Gpq{p}{q}/\Tpqf{p}{q}{f})$ and 
\begin{gather}
  (\Gpq{p}{q}/\Tpqf{p}{q}{f})/(\Gpqf{p}{q}{f}/\Tpqf{p}{q}{f}) 
  \cong \Gpq{p}{q}/\Gpqf{p}{q}{f}
  \label{eq:conj8}
\end{gather}
and the transversal of $\Tpqf{p}{q}{f}$ in $\Gpq{p}{q}$ spans $S$  over 
$\BR$ modulo~$f$.
\item[(x)] The stabilizer $\Gpqf{p}{q}{f}$ can be viewed as 
\begin{gather}
  \Gpqf{p}{q}{f} = \bigcap_{x \in \Tpqf{p}{q}{f}} C_{\Gpq{p}{q}}(x) 
                 = C_{\Gpq{p}{q}}(\Tpqf{p}{q}{f})
\end{gather}
where $C_{\Gpq{p}{q}}(x)$ is the centralizer of $x$ in $\Gpq{p}{q}$ and
$C_{\Gpq{p}{q}}(\Tpqf{p}{q}{f})$ is the centralizer of $\Tpqf{p}{q}{f}$ in
$\Gpq{p}{q}$.
\end{itemize}
\label{maintheorem}
\end{theorem}

\subsection{Summary of Some Basic Properties of Salingaros Vee Groups $G_{p,q}$}

In the following, we summarize some basic properties of  Salingaros vee groups $G_{p,q}$.
\begin{itemize}
\item $|G_{p,q}| =2^{1+p+q}$, $|G'_{p,q}|= 2$ because $G'_{p,q}= \{\pm 1\}$,
\item When $p+q\geq 1,$ $G_{p,q}$ is not simple as it has a nontrivial normal subgroup of order $2^m$ for every $m < 1+p+q$ (because every $p$-group of order $p^n$ has a normal subgroup of order $p^m$ for every $m \neq n$). 
\item When $p+q\geq 1,$ the center of any group $G_{p,q}$ is non-trivial since $2 \mid |Z(G_{p,q})|$ and so every group $G_{p,q}$ has a central element $\tau$ of order $2$. It is well-known that for any prime $p$ and a finite $p$-group $G \neq \{1\}$, the center of $G$ is non-trivial (Rotman~\cite{rotman}).
\item Every element of $G_{p,q}$ is of order $1,$ $2,$ or $4$.
\item Since $[G_{p,q}:G'_{p,q}] = |G_{p,q}|/|G'_{p,q}| = 2^{p+q},$ each $G_{p,q}$ 
has $2^{p+q}$ linear characters (James and Liebeck~\cite{james}).
\item The number $N$ of conjugacy classes in $G_{p,q}$, hence, the number of irreducible inequivalent representations of $G_{p,q}$, is $1+2^{p+q}$ (resp. $2+2^{p+q}$) when $p+q$ is even (resp. odd) 
(Maduranga~\cite{maduranga}).
\item We have the following result (see also Varlamov~\cite{varlamov}):
\begin{theorem}
\label{centerofGpq}
Let $G_{p,q} \subset \cl_{p,q}$. Then,
\begin{gather}
Z(G_{p,q})=
\begin{cases} 
\{\pm 1\} \cong \BZ_2 & \text{if $p-q \equiv 0,2,4,6 \pmod{8}$};\\
\{\pm 1,\pm \beta\} \cong \BZ_2 \times \BZ_2 & \text{if $p-q \equiv 1,5 \pmod{8}$};\\
\{\pm 1,\pm \beta\} \cong \BZ_4 & \text{if $p-q \equiv 3,7 \pmod{8}$}.
\end{cases}
\end{gather}
\label{ZGpqlemma}
\end{theorem}
\hspace*{-2ex}as a consequence of $Z(\cl_{p,q})= \{1\}$ (resp. $\{1,\beta\})$ when  $p+q$ is even resp. odd) where $\beta=\be_1\be_2\cdots \be_n,$ $n=p+q,$ is the unit pseudoscalar in~$\cl_{p,q}$.
\item In Salingaros' notation, the five isomorphism classes denoted as $N_{2k-1},N_{2k},\Omega_{2k-1},\Omega_{2k},S_k$ correspond to our notation $\Gpq{p}{q}$ as follows:
\begin{table}[h]
\begin{center}
\caption{Vee groups $G_{p,q}$ in Clifford algebras $\cl_{p,q}$}
\label{t1}
\renewcommand{\arraystretch}{1.0}
\begin{tabular}{ | c | c | c | c |} 
\hline
\multicolumn{1}{|c|}{Group}  & 
\multicolumn{1}{|c|}{Center} & 
\multicolumn{1}{|c|}{Group order}& 
\multicolumn{1}{|c|}{$\dim_\BR \cl_{p,q}$} 
\\ \hline
$N_{2k-1}$  
& \multicolumn{1}{|c|}{$\BZ_2$} 
& \multicolumn{1}{|c|}{$2^{2k+1}$}
& \multicolumn{1}{|c|}{$2^{2k}$} 
\\\hline
$N_{2k}$  
& \multicolumn{1}{|c|}{$\BZ_2$} 
& \multicolumn{1}{|c|}{$2^{2k+1}$}
& \multicolumn{1}{|c|}{$2^{2k}$} 
\\\hline
$\Omega_{2k-1}$  
& \multicolumn{1}{|c|}{$\BZ_2 \times \BZ_2$} 
& \multicolumn{1}{|c|}{$2^{2k+2}$}
& \multicolumn{1}{|c|}{$2^{2k+1}$} 
\\\hline
$\Omega_{2k}$  
& \multicolumn{1}{|c|}{$\BZ_2 \times \BZ_2$} 
& \multicolumn{1}{|c|}{$2^{2k+2}$}
& \multicolumn{1}{|c|}{$2^{2k+1}$} 
\\\hline
$S_{k}$  
& \multicolumn{1}{|c|}{$\BZ_4$} 
& \multicolumn{1}{|c|}{$2^{2k+2}$}
& \multicolumn{1}{|c|}{$2^{2k+1}$} 
\\\hline
\end{tabular}
\end{center}
\end{table}
\begin{align*}
N_{2k-1} & \leftrightarrow G_{p,q} \subset \cl_{p,q},\,\, p-q \equiv 0,2 \pmod 8,
\,\, \BK \cong \BR;\\[-0.5ex]
N_{2k}   & \leftrightarrow G_{p,q} \subset \cl_{p,q},\,\, p-q \equiv 4,6 \pmod 8,
\,\, \BK \cong \BH;\\[-0.5ex]
\Omega_{2k-1} & \leftrightarrow G_{p,q} \subset \cl_{p,q},\,\, p-q \equiv 1 \pmod 8,
\,\, \BK \cong \BR \oplus \BR;\\[-0.5ex]
\Omega_{2k} & \leftrightarrow G_{p,q} \subset \cl_{p,q},\,\, p-q \equiv 5 \pmod 8,
\,\, \BK \cong \BH \oplus \BH;\\[-0.5ex]
S_{k} & \leftrightarrow G_{p,q} \subset \cl_{p,q},\,\, p-q \equiv 3,7 \pmod 8,
\,\, \BK \cong \BC.
\end{align*}
(Salingaros~\cite{salingaros1,salingaros2,salingaros3}, Brown~\cite{brown}, Varlamov~\cite{varlamov})
\end{itemize}

\noindent
The first few vee groups $\Gpq{p}{q}$ of low orders  $4,8,16$ corresponding to Clifford algebras $\cl_{p,q}$ in dimensions $p+q=1,2,3$, are:
\begin{align*}
\mbox{Groups of order $4$:}\quad G_{1,0}&=D_4 ,\quad G_{0,1}=\BZ_4,\\
\mbox{Groups of order $8$:}\quad G_{2,0}&=D_8 = N_1,\quad G_{1,1}=D_8 = N_1,\quad G_{0,2}=Q_8 = N_2,\\
\mbox{Groups of order $16$:}\quad G_{3,0}&=S_1, \quad G_{2,1}=\Omega_1, \quad G_{1,2}=S_1, \quad G_{0,3}=\Omega_2.
\end{align*}
where $D_8$ is the dihedral group of a square, $Q_8$ is the quaternionic group, and 
$D_4 \cong \BZ_2 \times \BZ_2.$ For a construction of inequivalent irreducible representations and characters of these groups see~Maduranga and Ab\l amowicz~\cite{maduranga2}, and Maduranga~\cite{maduranga}.

\section{Central Product Structure of $G_{p,q}$}
We recall first a few definitions and results pertaining to finite $p$-groups that will be needed in the sequel.

\begin{definition}[Gorenstein \cite{gorenstein}]
A finite abelian $p$-group is \textit{elementary abelian} if every nontrivial element has order~$p$.
\end{definition}

\begin{example}($D_4=\BZ_2 \times \BZ_2$ is elementary abelian)\\
$(\BZ_p)^k = \BZ_p \times \cdots \times \BZ_p$ ($k$-times), in particular, $\BZ_2$, 
$\BZ_2 \times \BZ_2$, etc, are elementary abelian.
\label{ex:example3}
\end{example}

\begin{definition}[Dornhoff \cite{dornhoff}]
A finite $p$-group $P$ is \textit{extra-special} if (i) $P'=Z(P),$ (ii) $|P'| = p,$ and (iii) $P/P'$ is elementary abelian.
\end{definition}

\begin{example}($D_8$ is extra-special)\\
$D_8=\langle a,b \mid a^4=b^2=1, bab^{-1}=a^{-1}\rangle$ is extra-special because:
\begin{itemize}
\item $Z(D_8) = D_8'=[D_8,D_8] = \langle a^2 \rangle$, $|Z(D_8)|=2,$ 
\item $D_8/D_8'=D_8/Z(D_8) = 
\langle \langle a^2 \rangle, a\langle a^2 \rangle, 
b\langle a^2 \rangle, ab\langle a^2 \rangle\rangle \cong \BZ_2 \times \BZ_2.$
\end{itemize}
\label{ex:example4}
\end{example} 

\begin{example}($Q_8$ is extra-special)\\
$Q_8=\langle a,b \mid a^4=1, a^2=b^2, bab^{-1}=a^{-1}\rangle$ is extra-special because:
\begin{itemize}
\item $Z(Q_8) = Q_8'=[Q_8,Q_8] = \langle a^2 \rangle$, $|Z(Q_8)|=2,$ 
\item $Q_8/Q_8'=Q_8/Z(Q_8) = 
\langle \langle a^2 \rangle, a\langle a^2 \rangle, 
b\langle a^2 \rangle, ab\langle a^2 \rangle\rangle \cong \BZ_2 \times \BZ_2.$ 
\end{itemize}
\label{ex:example5}
\end{example}

Let us recall now definitions of internal and external central products of groups.
\begin{definition}[Gorenstein \cite{gorenstein}]
\leavevmode
\begin{enumerate}
\item A group $G$ is an \textit{internal central product} of two subgroups $H$ and $K$ if:
\begin{enumerate}
\item $[H,K] = \langle 1 \rangle$;
\item $G = HK$; 
\end{enumerate}
\item A group $G$ is an \textit{external central product} $H \circ K$ of two groups $H$ and $K$ with $H_{1} \leq Z(H)$ and $K_{1} \leq Z(K)$ if there exists an isomorphism 
$\theta :H_{1} \rightarrow K_{1}$ such that $G$ is $(H \times K)/N$ where 
$$
N = \lbrace (h,\theta (h^{-1})) \mid h \in H_{1} \rbrace.
$$ 
Clearly: $N \lhd (H \times K)$ and $|H \circ K| = |H||K|/|N| \leq |H \times K|=|H||K|.$
\end{enumerate}
\end{definition}

Here we recall an important result on extra-special $p$-groups as central products.
\begin{lemma}[Leedham-Green and McKay  \cite{mckay}]
An extra-special $p$-group has order $p^{2n + 1}$ for some positive integer~$n$, and is the iterated central product of non-abelian groups of order~$p^{3}$. 
\end{lemma}

As a consequence, we have the following lemma and a corollary. For their proofs, see~\cite{brown}.
\begin{lemma}
$Q_{8} \circ Q_{8} \cong D_{8} \circ D_{8} \ncong D_{8} \circ Q_{8}$, where $D_{8}$ is the dihedral group of order~8 and $Q_{8}$ is the quaternion group.
\label{lem:lem2}
\end{lemma}
\begin{corollary}
\leavevmode
\begin{itemize}
\item $\Gpq{3}{1} \cong D_8 \circ D_8 \cong Q_8 \circ Q_8,$
\item $\Gpq{4}{0} \cong D_8 \circ Q_8 \cong Q_8 \circ D_8.$
\end{itemize}
\end{corollary}

The following theorem is of critical importance for understanding the central product structure of Salingaros vee groups.

\begin{theorem}[Leedham-Green and McKay \cite{mckay}]
There are exactly two isomorphism classes of extra-special groups of order $2^{2n + 1}$ for positive integer $n$. One isomorphism type arises as the iterated central product of $n$ copies of $D_{8}$; the other as the iterated central product of $n$ groups isomorphic to $D_{8}$ and $Q_{8}$, including at least one copy of $Q_{8}$.
That is,
\begin{enumerate}
\item[1:] $D_{8} \circ D_{8} \circ \cdots \circ D_{8} \circ D_{8}$, or,
\item[2:] $D_{8} \circ D_{8} \circ \cdots \circ D_{8} \circ Q_{8}$.
\end{enumerate}
where it is understood that these are iterated central products; that is, $D_{8} \circ D_{8} \circ D_{8}$ is really $(D_{8} \circ D_{8}) \circ D_{8}$ and so on.
\label{ther:ther1}
\end{theorem}

Thus, the above theorem now explains the following theorem due to Salingaros regarding the iterative central product structure of the finite $2$-groups named after him.

\begin{theorem}[Salingaros Theorem \cite{salingaros3}]
Let $ {N_{1}} =  {D_{8}}$, 
    $ {N_{2}} =  {Q_{8}}$, and 
    $(G)^{\circ k}$ be the iterative central product 
    $G \circ G \circ \dots \circ G$ ($k$ times) of $G$. 
Then, for $k \geq 1$: 
\begin{enumerate}
\item $ {N_{2k-1}} \cong ( {N_1})^{\circ k} 
      =( {D_8})^{\circ k}$,
\item $ {N_{2k}} \cong ( {N_1})^{\circ k} \circ 
        {N_2} = ( {D_8})^{\circ (k-1)} \circ 
        {Q_8}$,
\item $\Omega_{2k-1} \cong  {N_{2k-1}} \circ 
( {\BZ_2 \times \BZ_2}) 
      =( {D_8})^{\circ k} \circ ( {\BZ_2 \times \BZ_2})$,
\item $\Omega_{2k} \cong  {N_{2k}} \circ 
( {\BZ_2 \times \BZ_2}) 
      =( {D_8})^{\circ (k-1)} \circ  {Q_8} 
       \circ ( {\BZ_2 \times \BZ_2})$,
\item $S_k \cong  {N_{2k-1}} \circ \BZ_4 \cong  {N_{2k}} \circ \BZ_4 
      =( {D_8})^{\circ k} \circ \BZ_4 \cong 
       ( {D_8})^{\circ (k-1)} \circ  {Q_8} \circ \BZ_4 $.
\end{enumerate}
\end{theorem}
\noindent
In the above theorem:
\begin{itemize}
\item $\BZ_2,$ $\BZ_4$ are cyclic groups of order~$2$ and~$4$, respectively; 
\item $D_8$ and $Q_8$ are the dihedral group of a square and the quaternionic group; 
\item $ {\BZ_2 \times \BZ_2}$ is elementary abelian of order~$4$;
\item $ {N_{2k-1}}$ and $ {N_{2k}}$ are extra-special of order $2^{2k+1}$;
e.g., $ {N_1}= {D_8}$ and $ {N_2}= {Q_8}$;
\item $\Omega_{2k-1},\Omega_{2k}, S_k$ are of order $2^{2k+2}$. 
\item $\circ$ denotes the iterative central product of groups with, e.g., 
      $({D_8})^{\circ k}$ denotes the iterative central product of $k$-copies of $ {D_8}$, etc.,
\end{itemize}
We can tabulate the above results for Salingaros vee groups $\Gpq{p}{q}$ of orders $\leq 256,$ $(p+q \leq 7)$ (Brown~\cite{brown}) in the following table:
\begin{table}[ht]
\renewcommand{\arraystretch}{1.5}
\caption{Salingaros Vee Groups $|\Gpq{p}{q}| \leq 256$}
\begin{center}
\begin{tabular}{|c|c|} \hline
Isomorphism Class & Salingaros Vee Groups  \\ \hline
$N_{2k}$ & $N_{0} \cong G_{0,0},\; N_{2} \cong Q_{8} \cong G_{0,2},\; N_{4} \cong G_{4,0},\; N_{6} \cong G_{6,0}$ \\ \hline
$N_{2k-1}$ & $N_{1} \cong D_{8} \cong G_{2,0},\; N_{3} \cong G_{3,1},\; N_{5} \cong G_{0,6}$ \\ \hline
$\Omega_{2k}$ & $\Omega_{0} \cong G_{1,0},\; \Omega_{2} \cong G_{0,3},\; \Omega_{4} \cong G_{5,0},\; \Omega_{6} \cong G_{6,1}$ \\ \hline
$\Omega_{2k-1}$ & $\Omega_{1} \cong G_{2,1},\; \Omega_{3} \cong G_{3,2},\; \Omega_{5} \cong G_{0,7}$ \\ \hline
$S_{k}$ & $S_{0} \cong G_{0,1},\; S_{1} \cong G_{3,0},\; S_{2} \cong G_{4,1},\; S_{3} \cong G_{7,0}$ \\ \hline
\end{tabular}
\end{center}
\end{table}

\section{Clifford Algebras Modeled with Walsh Functions}

Until now, the finite $2$-groups such as the Salingaros vee groups $\Gpq{p}{q}$ have appeared either as finite subgroups of the group of units $\clpq{p}{q}^{\times}$ in the Clifford algebra, or, as groups whose group algebra modulo a certain ideal generated by $1+\tau$ for some central element $\tau$ of order $2$ was isomorphic to the given Clifford algebra $\clpq{p}{q}.$ In these last two sections, we recall how the (elementary abelian) group $(\BZ_2)^n$ can be used to define a Clifford product on a suitable vector space.

In this section, we recall the well-known construction of the Clifford product on the set of monomial terms $\be_{\ua}$ indexed by binary $n$-tuples $\ua \in (\BZ_2)^n$, which, when extended by linearity, endows the set with the structure of the Clifford algebra $\clpq{p}{q}.$ This approach is described in 
Lounesto~\cite[Chapter 21]{lounesto}. We will show how it can be extended to Clifford algebras 
$\clpqr{p}{q}{r}$ over (real) quadratic vector spaces with degenerate quadratic forms. 

In the last section we will briefly mention the approach of Albuquerque and Majid~\cite{majid} in which the Clifford algebra structure is introduced in a suitably twisted group algebra $\BR^t[(\BZ_2)^n]$ using Hopf algebraic methods.

Let $\cb{B}^n=\{\underline{a} = a_1a_2\ldots a_n \mid a_i =0,1,\, \underline{a} \oplus \underline{b} = \underline{c} \mbox{ as } c_i=a_i+ b_i \bmod 2 \}$ be a group of binary $n$-tuples with addition $\oplus$, that is, $\cb{B}^n \cong (\BZ_2)^n$.

\begin{definition}[Walsh function] 
A \textit{Walsh function} $w_{\ua}$ indexed by $\ua \in \cb{B}^n$ is a function 
from $\cb{B}^n$ to the multiplicative group $\{\pm 1\}$ defined as
\begin{gather}
w_{\ua}(\ub) = (-1)^{\sum_{i=1}^{n} a_ib_i} = \pm 1, \quad \ua,\ub \in \cb{B}^n,
\end{gather}
which satisfies $w_{\uk}(\ua \oplus \ub) = w_{\uk}(\ua) w_{\uk}(\ub)$ and 
$w_{\ua}(\ub) = w_{\ub}(\ua)$. 
\end{definition}
Observe that the first condition on $w_{\ua}$ simply states that each $w_{\ua}$ is a group homomorphism from $\cb{B}^n$ to the group $\{\pm 1\}.$

\begin{definition}[Gray code]
A \textit{Gray code} $g: \cb{B}^n\rightarrow \cb{B}^n$ with the property 
$g(\ua \oplus \ub) = g(\ua) \oplus g(\ub)$ is defined as 
\begin{gather}
g(\uk)_1 = k_1, \quad g(\uk)_i=k_{i-1}+k_i \bmod 2, \quad i=2, \ldots, n.
\end{gather}
Thus, $g$ is a group isomorphism which reorders Walsh functions into a \textit{sequency order} with a 
\textit{single digit change code}~\cite[Section 21.2, page 281]{lounesto}.
\end{definition}
Given that the Gray code $g$ is an isomorphism, Lounesto defines its inverse $h: \cb{B}^n\rightarrow \cb{B}^n$ as
\begin{gather}
h(\ua)_i = \sum_{j=1}^i a_j \bmod 2.
\end{gather}
\noindent
Now, take an $\BR$-vector space $\cb{A}$ with a basis consisting of $2^n$ elements 
$\be_{\ua}$ labeled by the binary $n$-tuples $\ua=a_1a_2\ldots a_n$ as
\begin{gather}
\be_{\ua} = \be_1^{a_1}\be_2^{a_2} \cdots \be_n^{a_n}, \quad a_i = 0,1;
\end{gather}
for some $n$ symbols $\be_1,\be_2,\ldots,\be_n,$ and define an algebra product on 
$\cb{A}$ which on the basis elements $\be_{\ua}$ is computed as follows:
\begin{gather}
\be_{\ua}\be_{\ub} = (-1)^{\sum_{i=1}^p a_ib_i} w_{\ua}(h(\ub))\be_{\ua \oplus \ub},
\label{eq:eaeb} 
\end{gather}
for some $1\leq p \leq n.$ Then, together with this product, $\cb{A}$ becomes the Clifford algebra $\cl_{p,q}$, where $q=n-p$, of a non-degenerate quadratic form~$Q$ of signature $(p,q)$. See 
Lounesto \cite[Page 284]{lounesto} and his reference to~(\ref{eq:eaeb}) as the  formula of 
Brauer and Weyl from 1935~\cite{brauer}.

\begin{remark}
Observe that if the scalar factor in front of $\be_{\ua \oplus \ub}$ in~(\ref{eq:eaeb}) were set to be identically equal to~$1$, then we would have $\be_{\ua}\be_{\ub}=\be_{\ub}\be_{\ua}$ for any 
$\be_{\ua},\be_{\ub} \in \cb{A}.$ Thus, the algebra $\cb{A}$ would be isomorphic to the (abelian) group algebra $\BR[G]$ where $G \cong (\BZ_2)^n.$ That is, the scalar factor introduces a twist in the algebra product in $\cb{A}$ and so it makes~$\cb{A},$ hence the Clifford algebra $\clpq{p}{q},$ isomorphic to the twisted group algebra~$\BR^t[(\BZ_2)^n]$. 
\end{remark}

Formula~(\ref{eq:eaeb}) is encoded as a procedure \texttt{cmulWalsh3} in \texttt{CLIFFORD}, a Maple package for computations with Clifford algebras~\cite{clifford,parallelizing}. It has the following pseudo-code.

{\Fontv
\begin{lstlisting}[language = GAP]
cmulWalsh3:=proc(eI::clibasmon,eJ::clibasmon,B1::{matrix,list(nonnegint)}) 
local a,b,ab,monab,Bsig,flag,i,dim_V_loc,ploc,qloc,_BSIGNATUREloc; 
global dim_V,_BSIGNATURE,p,q;
options `Copyright (c) 2015-2016 by Rafal Ablamowicz and Bertfried Fauser. 
         All rights reserved.`;
if type(B1,list) then 
   ploc,qloc:=op(B1);
   dim_V_loc:=ploc+qloc:
   _BSIGNATUREloc:=[ploc,qloc]: 
   else 
   ploc,qloc:=p,q;    <<<-- this reads global p, q
   dim_V_loc:=dim_V: <<<-- this reads global dim_V
   _BSIGNATUREloc:=[ploc,qloc]:
   if not _BSIGNATURE=[ploc,qloc] then _BSIGNATURE:=[p,q] end if:
end if:
a:=convert(eI,clibasmon_to_binarytuple,dim_V_loc);
b:=convert(eJ,clibasmon_to_binarytuple,dim_V_loc);
ab:=oplus(a,b);
monab:=convert(ab,binarytuple_to_clibasmon);
return twist(a,b,_BSIGNATUREloc)*Walsh(a,hinversegGrayCode(b))*monab; 
end proc:
\end{lstlisting}
}

Now let us take a real quadratic vector space $(V,Q)$ with a degenerate quadratic form $Q$ such that 
$\dim V^{\perp} = r, $ while $Q$ restricted to the orthogonal complement of $V^{\perp}$ in $V$ has signature $(p,q)$, $(\dim V = n=p + q + r)),$ and we let a basis $\be_i, 1\leq i \leq n$ be such that 
$Q(\be_i)=1$ resp. $Q(\be_i)=-1$, resp. $Q(\be_i)=0$, for $0 \leq i \leq p,$ resp. $p+1 \leq i \leq p+q,$ resp. $p+q+1 \leq i \leq p+q+r$. We can now generate a universal Clifford algebra as the graded tensor product $\clpq{p}{q}{r} \cong \clpq{p}{q} \hotimes \bigwedge V^{\perp}$ with a Clifford product obtained by modifying the above formula~(\ref{eq:eaeb}) as follows: we introduce an extra scalar factor in front of 
$\be_{\ua \oplus \ub}$. This factor equals $1$ or, resp. $0$, depending whether the monomial elements 
$\be_{\ua}$ and $\be_{\ub}$ do not share, resp. do share, a common basis element $\be_i$ which squares to $0$ in $\clpqr{p}{q}{r}$, that is, such that $Q(\be_i)=0$. 

A modified pseudo-code of such procedure called \texttt{cmulWalsh3pqr} has been encoded in a new experimental package \texttt{eClifford} for computations in $\cl_{p,q,r}$~\cite{eclifford}.

{\Fontv
\begin{lstlisting}[language = GAP]
cmulWalsh3pqr:=proc(eI::eclibasmon,eJ::eclibasmon,B1::list(nonnegint)) 
local ploc,qloc,rloc,dim_V_loc,_BSIGNATUREloc,a,b,ab,monab,maxmaxindex,r_factor;
global twist,Walsh,hinverseGrayCode,oplus;
options `Copyright (c) 2015-2016 by Rafal Ablamowicz and Bertfried Fauser. 
         All rights reserved.`;
if nops(B1)=2 then  
   ploc,qloc:=op(B1);
   rloc:=0;
elif nops(B1)=3 then
   ploc,qloc,rloc:=op(B1);
else
   error `three non-negative integers p,q,r are needed in the list entered as
          the last argument but received \%1 instead`,B1
end if;
dim_V_loc:=ploc+qloc+rloc:
maxmaxindex:=max(op(eClifford:-eextract(eI)),op(eClifford:-eextract(eJ)));
if evalb(maxmaxindex>dim_V_loc) then
   error `maximum index \%1 found in the arguments \%2 and \%3 is larger 
          then dim_V = \%4 derived from the last argument \%5`,
          maxmaxindex,eI,eJ,dim_V_loc,B1
end if;
_BSIGNATUREloc:=[ploc,qloc]: 
a:=convert(eI,eclibasmon_to_binarytuple,dim_V_loc);
b:=convert(eJ,eclibasmon_to_binarytuple,dim_V_loc);
if rloc=0 then 
   r_factor:=1 
else
   r_factor:=mul((1+(-1)^(a[i]*b[i]))/2,i=ploc+qloc+1..(ploc+qloc+rloc)); 
end if;
if r_factor=0 then return 0 else
   ab:=oplus(a,b);
   monab:=convert(ab,binarytuple_to_eclibasmon);
   return twist(a,b,_BSIGNATUREloc)*Walsh(a,hinversegGrayCode(b))*monab; 
end if;
end proc:
\end{lstlisting}
}
\noindent
In the above, the code lines 25-33 accommodate the additional factor called \texttt{r\_factor} which equals $1$ or $0$ as indicated above\footnote{Note that such factor can also be computed by an \texttt{XOR} operation~\cite{private}.}. In particular, the Clifford algebra $\clpqr{0}{0}{n} \cong \bigwedge V$, the exterior (Grassmann) algebra.

\section{Clifford Algebras $\cl_{p,q}$ as Twisted Group Algebras}

In this last section we give a formal definition of a \textit{twisted group ring} (algebra) following 
Passman~\cite[Section 2]{passman}, and briefly refer to the paper by Albuquerque and Majid~\cite{majid} in which the authors discuss twisting of a real group algebra of $(\BZ_2)^n$ by using Hopf algebraic methods.



\begin{definition}[Passman \cite{passman}]
The \textit{twisted group ring} $k^t[G]$~\cite[Sect. 2]{passman} is an associative $k$-algebra, $k$ is a field, with a basis $\{\bar{x} \mid x \in G\}$ and multiplication defined distributively for all $x,y \in G$ as
\begin{gather} 
  \bar{x} \bar{y} = \gamma(x,y)\, \overline{xy}, \qquad \gamma(x,y) \in k^{\times} = k \setminus \{0\}.
\end{gather}
where the function $\gamma: G \times G \rightarrow k^{\times}$ satisfies
\begin{gather}
\gamma(x,y)\gamma(xy,z) = \gamma(y,z) \gamma(x,yz), \quad \forall z,y,z \in (\BZ_2)^n \quad \mbox{(cocycle condition)}
\end{gather}
to assure associativity 
$(\bar{x} \bar{y}) \, \bar{z} = \bar{x} \, (\bar{y} \, \bar{z})$ in $k^t[G]$ for any $x,y,z \in G.$
\end{definition}

\begin{lemma}[Passman \cite{passman}]
The following relations hold in $k^t[G]$.
\begin{itemize}
\item[(i)]  $\gamma(1,1)^{-1}\overline{1}$ is the identity in $k^t[G]$;
\item[(ii)] $\bar{x}^{-1}= \gamma(x,x^{-1})\gamma(1,1)^{-1} \overline{x^{-1}}
             = \gamma(x^{-1},x)\gamma(1,1)^{-1} \overline{x^{-1}},\, \forall x \in G$;
\item[(iii)] $(\bar{x}\bar{y})^{-1} = \bar{y}^{-1}\bar{x}^{-1},\, \forall x, y \in G$.
\end{itemize}
\label{lem:lemp}
\end{lemma}
\noindent
If $\gamma(1,1)=1$ in part (i) of the above lemma, then we call $\gamma$ \textit{normalized}, which can always be achieved by scaling. In part (ii), the inverse $\bar{x}^{-1}$ is the result of the action of the antipode on $\bar{x}$ in the Hopf algebra sense, or, it can be viewed as the (un-normalized) action of the transposition map $\tp$ introduced in \cite{ablamowicz1,ablamowicz2,ablamowicz3} and mentioned in Section~\ref{sub:sectt}.

For a Hopf algebraic discussion of Clifford algebras $\cl_{p,q}$ as twisted group
algebras $\BR^t[(\BZ_2)^n]$, where the twisting is accomplished via a $2$-cocycle $F$ which twists the group algebra $k[(\BZ_2)^n]$ into a cotriangular Hopf algebra with a suitable cotriangular 
structure~$\cb{R}$, see~\cite{albuquerque,downs} and references therein. Note that if $\gamma$ is trivial, then the twist is trivial and the twisted group algebra is just the group algebra $k[G]$; 
if it is given by the \texttt{XOR} function on binary tuples, we get the Grassmann product (including a graded tensor product, or a graded switch; if $\gamma$ is the choice described by Lounesto 
in~(\ref{eq:eaeb}), we get the Clifford algebra $\clpq{p}{q}$~\cite{private}.

\section{Conclusions}
As stated in the Introduction, the main goal of this survey paper has been to collect and summarize properties of certain finite $2$-groups which appear in Clifford algebras~$\clpq{p}{q}$. On one hand, these Salingaros-defined groups $\Gpq{p}{q}$ appear as subgroups of the group of invertible elements. These subgroups play an important role in relation to the set of orthogonal primitive idempotents, with the help of which one defines spinorial representations. It has been observed by Salingaros, that these groups belong to five non-isomorphic families. On the other hand, one knows that all Clifford algebras $\clpq{p}{q}$ are classified into five different families of simple and semisimple algebras depending on the values of $(p,q)$ and $p+q$ (the Periodicity of Eight). Another connection with finite Salingaros groups appears via Chernov's observation that the algebras $\clpq{p}{q}$ can be viewed as images of group algebras, most likely of the groups $\Gpq{p}{q}$ modulo a suitable ideal generated by a central nontrivial idempotent in the group algebra. This shows that the theory of extra-special $2$-groups has a direct bearing on the structure of the Clifford algebras $\clpq{p}{q}$. Finally, we have observed how Clifford algebras can be obtained by twisting a group algebra of $(\BZ_2)^n$, either by using the Walsh functions, or equivalently but in a more sound mathematical way, by using a $2$-cocycle and the formalism of cotriangular Hopf algebras \cite{downs}.

\section{Acknowledgments}
Author of this paper is grateful to Dr. habil. Bertfried Fauser for his remarks and comments which have 
helped improve this paper.

\appendix
\section{The Structure Theorem on Clifford Algebras}

In this appendix we list the main structure theorem for real Clifford algebras $\clpq{p}{q}.$ For more information on Clifford algebras, see~\cite{chevalley,lam,lounesto}.

\begin{StructureTheorem}
Let $\cl_{p,q}$ be the universal Clifford algebra over $(V,Q)$, $Q$ is non-degenerate of signature $(p,q)$.
\begin{itemize}
\item[(a)] When $p-q \neq 1 \bmod 4$ then $\cl_{p,q}$ is a simple algebra of dimension $2^{p+q}$ isomorphic with a full matrix algebra $\Mat(2^k, \BK)$ over a division ring $\BK$ where $k = q - r_{q-p}$ and $r_i$ is the Radon-Hurwitz number.\footnote{The Radon-Hurwitz number is defined by recursion as $r_{i+8}=r_i+4$ and these initial values: $r_0=0,$ $r_1=1,$ $r_2=r_3=2,r_4=r_5=r_6=r_7=3.$ \label{page:pageRH}} Here $\BK$ is one of $\BR, \BC$ or $\BH$ when $(p-q) \bmod 8$ is $0,2,$ or $3,7$, or $4,6$.
\item[(b)] When $p-q = 1 \bmod 4$ then $\cl_{p,q}$ is a semisimple algebra of dimension $2^{p+q}$ isomorphic to $\Mat(2^{k-1}, \BK) \oplus \Mat(2^{k-1}, \BK),$ $k = q - r_{q-p}$, and $\BK$ is isomorphic to $\BR$ or $\BH$ depending whether $(p-q) \bmod 8$ is $1$ or~$5$. Each of the two simple direct components of $\cl_{p,q}$ is projected out by one of the two central idempotents $\frac12(1\pm \be_{12 \ldots n}).$
\item[(c)] Any element $f$ in $\cl_{p,q}$ expressible as a product
\begin{equation}
f = \frac12(1\pm \be_{\iu_1})\frac12(1\pm \be_{\iu_2})\cdots\frac12(1\pm \be_{\iu_k})
\label{eq:f}
\end{equation}
where $\be_{\iu_j},$ $j=1,\ldots,k,$ are commuting basis monomials in $\cb{B}$ with square $1$ and $k = q - r_{q-p}$ generating a group of order $2^k$, is a primitive idempotent in $\cl_{p,q}.$ Furthermore, $\cl_{p,q}$ has a complete set of $2^k$ such primitive mutually annihilating idempotents which add up to the unity $1$ of $\cl_{p,q}$. 
\item[(d)] When $(p-q) \bmod 8$ is $0,1,2,$ or $3,7$, or $4,5,6$, then the division ring $\BK = f \cl_{p,q}f$ is isomorphic to $\BR$ or $\BC$ or $\BH$, and the map $S \times \BK \rightarrow S,$ $(\psi,\lambda) \mapsto \psi\lambda$ defines a right  $\BK$-module structure on the minimal left ideal $S=\cl_{p,q}f.$ 
\item[(e)] When $\cl_{p,q}$ is simple, then the map
\begin{equation}
\cl_{p,q} \stackrel{\gamma}{\longrightarrow} \End_\BK(S), \quad u \mapsto \gamma(u), \quad \gamma(u)\psi = u \psi
\label{eq:simple}
\end{equation}
gives an irreducible and faithful representation of $\cl_{p,q}$ in $S.$
\item[(f)] When $\cl_{p,q}$ is semisimple, then the map
\begin{equation}
\cl_{p,q} \stackrel{\gamma}{\longrightarrow} \End_{\BK \oplus \hat{\BK}}(S \oplus \hat{S}), \quad u \mapsto \gamma(u), \quad \gamma(u)\psi = u \psi
\label{eq:semisimple}
\end{equation}
gives a faithful but reducible representation of $\cl_{p,q}$ in the double spinor space $S \oplus \hat{S}$ where $S =\{ u f\, |\, u \in \cl_{p,q}\}$, $\hat{S} =\{ u \hat{f}\, |\, u \in \cl_{p,q}\}$ and 
$\hat{\phantom{m}}$ stands for the grade-involution in $\cl_{p,q}.$ In this case, the ideal $S \oplus \hat{S}$ is a right $\BK \oplus \hat{\BK}$-module structure,  
$\hat{\BK} = \{\hat{\lambda} \,| \, \lambda \in \BK \}$, and $\BK \oplus \hat{\BK}$ is isomorphic to $\BR \oplus \BR$ when $p-q=1 \bmod 8$ or to $\BH \oplus \hat{\BH} $ when 
$p-q=5 \bmod 8.$
\end{itemize}
\label{th:structure}
\end{StructureTheorem}

\end{document}